\documentclass[]{scrartcl}

\usepackage[]{amsmath, amssymb,amsfonts,amsthm,mathtools,braket,color,enumerate}
\usepackage{mathrsfs}
\usepackage{booktabs}

\usepackage[margin=25mm]{geometry}

\usepackage{hhline}
\usepackage{url}
\usepackage{here}
\usepackage{stmaryrd}
\SetSymbolFont{stmry}{bold}{U}{stmry}{m}{n}

\usepackage{zref-clever}
\zcsetup{cap=true}
\AddToHook{env/proposition/begin}{\zcsetup{countertype={theorem=proposition}}}
\AddToHook{env/corollary/begin}{\zcsetup{countertype={theorem=corollary}}}
\AddToHook{env/lemma/begin}{\zcsetup{countertype={theorem=lemma}}}
\AddToHook{env/example/begin}{\zcsetup{countertype={theorem=example}}}
\AddToHook{env/conjecture/begin}{\zcsetup{countertype={theorem=conjecture}}}
\zcRefTypeSetup{conjecture}{
  Name-sg = Conjecture, name-sg = conjecture,
  Name-pl = Conjectures, name-pl = conjectures
}

\usepackage{hyperref}
\hypersetup{
  colorlinks   = true, 
  urlcolor     = blue, 
  linkcolor    = blue, 
  citecolor   = red 
}

\usepackage{tikz-cd}
\usepackage{tikz}
\usetikzlibrary{decorations}
\usetikzlibrary{arrows}
\tikzstyle{v} = [circle, draw, inner sep=2pt, minimum size=3pt, fill=black]
\tikzstyle{l} = [rectangle, draw, rounded corners]
\tikzstyle{c} = [circle, draw, inner sep=1pt, minimum size=14pt]

\theoremstyle{plain}
\newtheorem{theorem}{Theorem}[section]

\theoremstyle{definition}

\newtheorem{conjecture}[theorem]{Conjecture}

\DeclareMathOperator{\Kernel}{Ker}

\title{A counterexample to Terao's conjecture for free hyperplane arrangements over  $\mathbb{F}_{121}$}
\author{
Shuhei Tsujie
\thanks{Department of Mathematics, Hokkaido University of Education, Asahikawa, Hokkaido 070-8621, Japan. \\
E-mail: tsujie.shuhei@a.hokkyodai.ac.jp}
}

\date{}

\begin{document}
\maketitle

\begin{abstract}
We give a counterexample to Terao's conjecture for free hyperplane arrangements defined over the finite field $\mathbb{F}_{121}$. 
\end{abstract}




\section{Introduction}
Let $V$ be an $\ell$-dimensional vector space over a field $\mathbb{K}$. 
A central hyperplane arrangement $\mathcal{A}$ is a finite collection of linear hyperplanes in $V$. 
Let $S = \mathbb{K}[x_{1}, \dots, x_{\ell}]$ be the coordinate ring of $V$. 
A logarithmic derivation module $D(\mathcal{A})$ is the set of $\mathbb{K}$-derivations $\theta \colon S \to S$ such that $\theta(\alpha_{H}) \in \alpha_{H}S$ for all $H \in \mathcal{A}$, where $\alpha_{H}$ denotes a fixed linear form such that $\Kernel \alpha_{H} = H$. 
We say that $\mathcal{A}$ is free with exponents $(d_{1}, \dots, d_{\ell})$ if $D(\mathcal{A})$ is a free $S$-module with a homogeneous basis whose degrees are $(d_{1}, \dots, d_{\ell})$. 
Let $L(\mathcal{A})$ denote the intersection lattice of $\mathcal{A}$, that is, the set of all intersections of hyperplanes in $\mathcal{A}$ ordered by reverse inclusion. 

\begin{conjecture}[Terao {\cite[Problem 1]{terao1983exponents-spc}}, {\cite[Conjecture 4.138]{orlik1992arrangements}}]\label{Terao's conjecture}
Let $\mathcal{A}$ and $\mathcal{B}$ be arrangements over a fixed field $\mathbb{K}$. 
Suppose that $\mathcal{A}$ is free and the intersection lattices $L(\mathcal{A})$ and $L(\mathcal{B})$ are isomorphic. 
Then $\mathcal{B}$ is free. 
\end{conjecture}

For arrangements of rank $3$ with at most $14$ hyperplanes, this conjecture was verified in arbitrary characteristic by Barakat and K{\"u}hne \cite{barakat2023computing-moc}. 
Note that the condition that two arrangements $\mathcal{A}$ and $\mathcal{B}$ are defined over the same field is necessary. 
In fact, Ziegler \cite[Example 4.1]{ziegler1990matroid-totams} constructed the following example: 
\begin{itemize}
\item $\mathcal{A}$ has a defining polynomial with coefficients in the ring of Eisenstein integers $\mathbb{Z}[\omega]$, where $\omega$ is a primitive cube root of unity, and is free with exponents $(1,4,4)$. 
\item Let $\mathfrak{p} = (1-\omega)$ be the prime ideal of $\mathbb{Z}[\omega]$ dividing $3$. 
Let $\mathcal{A}_{\mathfrak{p}}$ be the arrangement over $\mathbb{F}_{3} \simeq \mathbb{Z}[\omega]/\mathfrak{p}$ obtained by reducing the defining form of each hyperplane in $\mathcal{A}$ modulo $\mathfrak{p}$. 
The arrangement $\mathcal{A}_{\mathfrak{p}}$ is non-free. 
\item $L(\mathcal{A})$ and $L(\mathcal{A}_{\mathfrak{p}})$ are isomorphic. 
\end{itemize}

It is a key point that the Frobenius derivation of degree $3$ belongs to $D(\mathcal{\mathcal{A}_{\mathfrak{p}}})$ since $\mathcal{\mathcal{A}_{\mathfrak{p}}}$ is defined over $\mathbb{F}_{3}$. 
If $\mathcal{\mathcal{A}_{\mathfrak{p}}}$ were free, Terao's factorization theorem would imply that the exponents of $\mathcal{\mathcal{A}_{\mathfrak{p}}}$ are $(1,4,4)$. 
Every homogeneous derivation of degree $3$ in $D(\mathcal{A}_{\mathfrak{p}})$ would be a multiple of the Euler derivation. 
However, the Frobenius derivation of degree $3$ is not a multiple of the Euler derivation. 
This contradiction shows that $\mathcal{\mathcal{A}_{\mathfrak{p}}}$ is non-free. 

In this article, we will construct the following arrangements: 
\begin{itemize}
\item $\mathcal{A}$ has a defining polynomial with coefficients in the ring $\mathbb{Z}[\psi]$, where $\psi$ is the supergolden ratio, and $\mathcal{A}$ is free with exponents $(1,13,13)$. 
\item In $\mathbb{Z}[\psi]$, the ideal $(11)$ factors as $(11) = \mathfrak{p}_{1}\mathfrak{p}_{2}$. 
Let $\mathcal{A}_{\mathfrak{p}_{i}}$ be the arrangement obtained by reducing the defining form of each hyperplane in $\mathcal{A}$ modulo $\mathfrak{p}_{i}$ for $i \in \{1,2\}$. 
\item $\mathcal{A}_{\mathfrak{p}_{1}}$ is defined over $\mathbb{F}_{11} \simeq \mathbb{Z}[\psi]/\mathfrak{p}_{1}$ and non-free. 
We can regard $\mathcal{A}_{\mathfrak{p}_{1}}$ as an arrangement
over $\mathbb{F}_{121}$. 
\item $\mathcal{A}_{\mathfrak{p}_{2}}$ is defined over $\mathbb{F}_{121} \simeq \mathbb{Z}[\psi]/\mathfrak{p}_{2}$ and free.  
\item $L(\mathcal{A}), L(\mathcal{A}_{\mathfrak{p}_{1}})$, and $L(\mathcal{A}_{\mathfrak{p}_{2}})$ are isomorphic. 
\end{itemize}

In particular, we have the following theorem. 
\begin{theorem}
The arrangements $\mathcal{A}_{\mathfrak{p}_{1}}$ and $\mathcal{A}_{\mathfrak{p}_{2}}$ give a counterexample to \zcref{Terao's conjecture} in the case $\mathbb{K} = \mathbb{F}_{121}$. 
\end{theorem}

\section{Construction}

Let $\psi$ be the supergolden ratio, that is, the real root of the cubic polynomial $x^{3}-x^{2}-1$. 
Define the $3$-dimensional arrangement $\mathcal{A}$ over $\mathbb{Z}[\psi]$ by 

\begin{align*}
\begin{aligned}
Q(\mathcal{A}) ={}& xyz \\
&{}\times (x-y)(y-z)(-\psi x+z) \\
&{}\times (x-z)(-\psi x+y)(-\psi y+z) \\
&{}\times (x-y+z)(\psi x+y-z)(-\psi x+\psi y+z) \\
&{}\times (x-\psi^{2}y+z)(\psi x+y-\psi^{2}z)(-\psi^{3}x+\psi y+z) \\
&{}\times (x-\psi^{2}y+\psi^{2}z)(\psi^{3}x+y-\psi^{2}z)(-\psi^{3}x+\psi^{3}y+z) \\
&{}\times (x+\psi y-\psi z)(-\psi^{2}x+y+\psi z)(\psi^{2}x-\psi^{2}y+z) \\
&{}\times (x+\psi y-\psi^{3}z)(-\psi^{4}x+y+\psi z)(\psi^{2}x-\psi^{4}y+z) \\
&{}\times (x+\psi^{3}y-\psi^{3}z)(-\psi^{4}x+y+\psi^{3}z)(\psi^{4}x-\psi^{4}y+z).
\end{aligned}
\end{align*}
Using SageMath \cite{sagemath}, we can check that $\mathcal{A}$ is free with exponents $(1,13,13)$. 

In the ring $\mathbb{Z}[\psi]$, the ideal $(11)$ factors as 
\begin{align*}
(11) = \mathfrak{p}_{1} \mathfrak{p}_{2}, \text{ where } 
\mathfrak{p}_{1} = (11, \psi-5) \text{ and }
\mathfrak{p}_{2} = (11, \psi^{2}+4\psi+9). 
\end{align*}

Therefore, the defining polynomial $Q(\mathcal{A}_{\mathfrak{p}_{1}}) \in \mathbb{F}_{11}[x,y,z]$ of $\mathcal{A}_{\mathfrak{p}_{1}}$ is 
\begin{align*}
\begin{aligned}
Q(\mathcal{A}_{\mathfrak{p}_{1}}) ={}& xyz \\
&{}\times (x-y)(y-z)(-5 x+z) \\
&{}\times (x-z)(-5 x+y)(-5 y+z) \\
&{}\times (x-y+z)(5 x+y-z)(-5 x+5 y+z) \\
&{}\times (x-5^{2}y+z)(5 x+y-5^{2}z)(-5^{3}x+5 y+z) \\
&{}\times (x-5^{2}y+5^{2}z)(5^{3}x+y-5^{2}z)(-5^{3}x+5^{3}y+z) \\
&{}\times (x+5 y-5 z)(-5^{2}x+y+5 z)(5^{2}x-5^{2}y+z) \\
&{}\times (x+5 y-5^{3}z)(-5^{4}x+y+5 z)(5^{2}x-5^{4}y+z) \\
&{}\times (x+5^{3}y-5^{3}z)(-5^{4}x+y+5^{3}z)(5^{4}x-5^{4}y+z).
\end{aligned}
\end{align*}

Moreover, the defining polynomial $Q(\mathcal{A}_{\mathfrak{p}_{2}}) \in \mathbb{F}_{121}[x,y,z]$ of $\mathcal{A}_{\mathfrak{p}_{2}}$ is 
\begin{align*}
\begin{aligned}
Q(\mathcal{A}_{\mathfrak{p}_{2}}) ={}& xyz \\
&{}\times (x-y)(y-z)(-s x+z) \\
&{}\times (x-z)(-s x+y)(-s y+z) \\
&{}\times (x-y+z)(s x+y-z)(-s x+s y+z) \\
&{}\times (x-s^{2}y+z)(s x+y-s^{2}z)(-s^{3}x+s y+z) \\
&{}\times (x-s^{2}y+s^{2}z)(s^{3}x+y-s^{2}z)(-s^{3}x+s^{3}y+z) \\
&{}\times (x+s y-s z)(-s^{2}x+y+s z)(s^{2}x-s^{2}y+z) \\
&{}\times (x+s y-s^{3}z)(-s^{4}x+y+s z)(s^{2}x-s^{4}y+z) \\
&{}\times (x+s^{3}y-s^{3}z)(-s^{4}x+y+s^{3}z)(s^{4}x-s^{4}y+z), 
\end{aligned}
\end{align*}
where $s$ is an element satisfying $\mathbb{F}_{121} = \mathbb{F}_{11}(s)$ and $s^{2}+4s+9=0$. 
Using SageMath \cite{sagemath}, we can verify that $\mathcal{A}_{\mathfrak{p}_{1}}$ is non-free, $\mathcal{A}_{\mathfrak{p}_{2}}$ is free, and the three lattices $L(\mathcal{A}), L(\mathcal{A}_{\mathfrak{p}_{1}})$, and $L(\mathcal{A}_{\mathfrak{p}_{2}})$ are isomorphic.  

 \section*{Acknowledgments}
S. T. was supported by JSPS KAKENHI Grant Numbers JP23H00081 and JP26K16955. 

\section*{AI Disclosure}
The author used ChatGPT (GPT-6 Astra) to search for hyperplane arrangements.
The author independently verified the results and takes responsibility for their correctness.

\bibliographystyle{amsplain}
\bibliography{bibfile}

\end{document}